\documentclass{article}
\usepackage{spconf}
\usepackage[named]{algo}
\usepackage{graphics}
\usepackage{graphicx}
\usepackage{subfigure}
\usepackage{stfloats}
\usepackage{amsmath}
\usepackage{amsfonts}

\newcommand{\btheta}{\mbox{\boldmath $\theta$}}

\newcommand{\bbeta}{\mbox{\boldmath  $\beta$}}

\newcommand{\balpha}{\mbox{\boldmath $\alpha$}}

\newcommand{\bPsi}{\mbox{\boldmath   $\Psi$}}

\newcommand{\blambda}{\mbox{\boldmath $\lambda$}}

\ninept

\begin{document}

\title{Fast Frame-Based Image Deconvolution Using \\ Variable
Splitting and Constrained Optimization}

\name{M\'{a}rio A. T. Figueiredo,  Jos\'{e} M. Bioucas-Dias, and Manya V. Afonso}
\address{ Instituto de Telecomunica\c{c}\~{o}es, \\
Instituto Superior T\'{e}cnico, Technical University of Lisbon, Portugal\\
Email: $\{$mario.figueiredo,\ jose.bioucas, \ mafonso$\}$@lx.it.pt\thanks{M. Afonso
is supported by a EU Marie-Curie Fellowship (EST-SIGNAL program: {\scriptsize\tt est-signal.i3s.unice.fr});
contract MEST-CT-2005-021175.}}

\maketitle

\begin{abstract}
We propose a new fast algorithm for solving one of the standard formulations of frame-based image deconvolution: an unconstrained optimization problem, involving an $\ell_2$ data-fidelity term and a non-smooth regularizer. Our approach is based on using variable splitting to obtain an equivalent constrained optimization formulation, which is then
addressed with an augmented Lagrangian method. The resulting algorithm efficiently uses a regularized version of the Hessian of the data fidelity term, thus exploits second order
information. Experiments on a set of image deblurring benchmark problems show that our
algorithm is clearly faster than previous state-of-the-art methods.
\end{abstract}

\section{Introduction} \subsection{Problem Formulation}
\label{sec:intro} The standard model in image deblurring
assumes that the noisy  blurred observed version ${\bf y}$,
of an original image ${\bf x}\in\mathbb{R}^n$, was obtained via
\[
{\bf y} = {\bf H}{\bf x} + {\bf w},
\]
where ${\bf H}$  is the matrix representation of a convolution
and ${\bf w}$ is Gaussian white noise. In frame-based
deblurring/deconvolution, the unknown image ${\bf x}$ is expressed
as ${\bf x} = {\bf W}\bbeta$, where the columns of matrix ${\bf W}$
are the elements of a frame, such as a wavelet orthonormal basis
or a redundant dictionary \cite{ChauxCombettesPesquetWajs}, \cite{CombettesSIAM},
\cite{DaubechiesDefriseDeMol},  \cite{EladCVPR2006},
\cite{FigueiredoNowak2003}, \cite{FigueiredoNowak2005}. The coefficients
of this representation are then estimated, under one of the
well-known sparsity inducing regularizers, typically the $\ell_1$ norm,
leading to the optimization problem
\begin{equation}
\widehat{\bbeta} = \arg\min_{\bbeta} \frac{1}{2}\| {\bf HW}\bbeta -
{\bf y}\|_2^2 + \tau \, \phi (\bbeta); \label{problem}
\end{equation}
in (\ref{problem}), $\phi:\mathbb{R}^m \! \rightarrow \overline{\mathbb{R}}$ is
the {\it regularizer}, which is usually convex but nonsmooth,
and $\tau \geq 0$ is the regularization parameter \cite{CombettesSIAM}.
This formulation is called the {\it synthesis approach}
\cite{EladMilanfarRubinstein}, since it is based on the synthesis
equation ${\bf x} = {\bf W}\bbeta$. In the last decade, a considerable
amount of research has been devoted  to designing efficient algorithms
for solving (\ref{problem}). This interest has been further stimulated
by the recent emergence of compressive sensing (CS) \cite{Candes},
\cite{donoho}, since CS reconstruction can be formulated as
(\ref{problem}) \cite{Haupt}, \cite{Zhu}.

\subsection{Previous Algorithms} In most practical problems (including CS),
matrix ${\bf HW}$ (and even ${\bf H}$ or ${\bf W}$) cannot be stored
explicitly and it is highly impractical to access portions (lines/columns or
blocks) of it.  These facts preclude most off-the-shelf optimization
algorithms from being directly used and has stimulated the
development of special purpose methods. These methods operate under the constraint
${\bf H}$ and ${\bf W}$ (and their transposes) can only be used to form
matrix-vector products, since these products can be
performed efficiently using the FFT and fast wavelet transforms.

Arguably, the standard algorithm for solving (\ref{problem}) is the
so-called {\it iterative shrinkage/thresholding} (IST), which can be
derived from different viewpoints: expectation-maximization
\cite{FigueiredoNowak2003}, {\it majorization-minimization}
\cite{DaubechiesDefriseDeMol}, \cite{FigueiredoNowak2005},
forward-backward operator splitting \cite{CombettesSIAM}, \cite{Hale}.
A key ingredient of IST is the so-called
shrinkage/thresholding function associated to $\phi$, $\bPsi_{\tau\phi}:\mathbb{R}^n\rightarrow \mathbb{R}^n$,
defined as
\begin{equation}
\bPsi_{\tau\phi}({\balpha}) = \arg\min_{\bbeta} \frac{1}{2}\|\balpha-\bbeta\|_2^2 + \tau\phi(\bbeta).
\label{MPM}
\end{equation}
An excellent coverage of these functions, also known as {\it
Moreau proximal maps}, can be found in \cite{CombettesSIAM}.

The fact that IST tends to be slow, in particular when ${\bf H}$
is poorly conditioned, has stimulated some recent
research aimed at obtaining faster variants. The recent
{\it two-step IST} (TwIST) algorithm \cite{TwIST}, in which each
iteration uses the two previous iterates (rather than only
the previous one, as in IST), was shown to be considerably faster
than IST on various deconvolution problems. Another
two-step variant of IST, named {\it fast IST algorithm} (FISTA),
was recently proposed and also shown to be faster than IST \cite{FISTA}.
A recent strategy to obtaining faster variants of IST
consists in using more aggressive choices of step size in each
iteration. This is the case in the SpaRSA ({\it sparse reconstruction
by separable approximation}) framework \cite{SpaRSA_ICASSP},
\cite{SpaRSA_SP}, which was also shown to clearly outperform
standard IST.

\subsection{Proposed Approach}
The approach proposed in this paper is based on
variable splitting. The idea is to split the variable $\bbeta$ into
a pair of variables $\bbeta$ and $\btheta$, each to serve as the
argument of each of the two functions in (\ref{problem}), and
then minimize the sum of the two functions under the constraint
that the two variables have to be equal, thus making the problems
equivalent. This rationale has been  recently used in the
split-Bregman method \cite{GoldsteinOsher}, which was proposed to
address constrained optimization formulations for solving
inverse problems. In this paper, we exploit a different
splitting to attack problem (\ref{problem}), arguably the most
classical formulation for frame-based regularization of
linear inverse problems \cite{ChauxCombettesPesquetWajs},
\cite{CombettesSIAM}.

The constrained optimization problem produced by the
splitting procedure is addressed using an augmented
Lagrangian (AL) algorithm \cite{NocedalWright}.
AL was shown to be equivalent to the Bregman iterative
methods \cite{Setzer}, \cite{YinOsherGoldfarbDarbon}. We adopt
the AL perspective, rather than the Bregman  view, as it is a
more standard optimization tool. We show that by exploiting the
fact that ${\bf W}$ is a frame, the resulting algorithm solves
(\ref{problem}) much faster than the previous state-of-the-art
methods FISTA \cite{FISTA}, TwIST \cite{TwIST}, and SpaRSA \cite{SpaRSA_SP}.

The speed of the proposed algorithm may be justified by the
fact that it uses (a regularized version of) the Hessian of the data
fidelity term, ${\bf W}^T{\bf H}^T{\bf H\, W}$,
while the above mentioned algorithms essentially only use gradient
information. Although, as referred earlier, this matrix can not be
formed, we show that if ${\bf W}$ is a tight frame and ${\bf H}$ a
convolution, our algorithm can use it in an efficient way.

\section{Basic Tools}

\subsection{Variable Splitting} Consider an unconstrained
optimization problem in which the objective is the sum of
two functions:
\begin{equation}
\min_{{\bf u}\in \mathbb{R}^n} f_1({\bf u}) + f_2({\bf u}).\label{unconstrained_basic}
\end{equation}
Variable splitting (VS) is a simple procedure in which a
new variable ${\bf v}$ is introduced to serve as the argument of
$f_2$, under the constraint that ${\bf u} = {\bf v}$. In other
words, the constrained problem
\begin{equation}\begin{array}{cl}
{\displaystyle \min_{{\bf u},{\bf v}\in \mathbb{R}^n}} & f_1({\bf u}) + f_2({\bf v})\\
\mbox{s.t.} & {\bf u} = {\bf v},
\end{array}\label{constrained_basic}
\end{equation}
is equivalent to  (\ref{unconstrained_basic}), since in
the feasible set $\{({\bf u},{\bf v}): {\bf u} = {\bf v}\}$,
the objective function in (\ref{constrained_basic}) coincides with that in
(\ref{unconstrained_basic}).

VS was used in \cite{Wang} to derive a fast
algorithm for total-variation based restoration. VS was
also used in \cite{BioucasFigueiredo2008} to handle problems  where
instead of the single regularizer $\tau \phi(\bbeta)$ in (\ref{problem}),
there is a linear combination of two (or more) regularizers: $\tau_1
\phi_1(\bbeta) + \tau_2 \phi_2(\bbeta)$. In
\cite{BioucasFigueiredo2008} and \cite{Wang}, the constrained
problem (\ref{constrained_basic}) is attacked by a quadratic penalty
approach, i.e., by solving
\begin{equation}
 \min_{{\bf u},{\bf v}\in \mathbb{R}^n}  f_1({\bf u}) + f_2({\bf v}) +
 \frac{\mu}{2}\, \|{\bf u} -{\bf v}\|_2^2,
\label{quadratic_penalty}
\end{equation}
by alternating minimization with respect to ${\bf u}$ and ${\bf v}$,
while slowly increasing $\mu$ to force the solution of
(\ref{quadratic_penalty}) to approach that of
(\ref{constrained_basic}). The idea is that each step
of this alternating minimization may be much easier than
the original unconstrained problem (\ref{unconstrained_basic}).
The drawback is that as $\mu$ increases, the intermediate
minimization problems become increasingly ill-conditioned,
thus causing numerical problems \cite{NocedalWright}.

A similar VS approach underlies the recently proposed
split-Bregman methods \cite{GoldsteinOsher}. In those
methods, the constrained problem (\ref{constrained_basic}) is
addressed using a Bregman iterative algorithm, which has been shown to
be equivalent to the AL method \cite{YinOsherGoldfarbDarbon}.

\subsection{Augmented Lagrangian}
Consider a linear equality constrained optimization problem
\begin{equation}\begin{array}{cl}
 {\displaystyle \min_{{\bf z}\in \mathbb{R}^d}} & E({\bf z})\\
 \mbox{s.t.} & {\bf A z - b = }\mbox{\boldmath $0$}, \end{array}\label{constrained_linear}
\end{equation}
where ${\bf b} \in \mathbb{R}^p$ and ${\bf A}\in \mathbb{R}^{p\times
d}$. The so-called augmented Lagrangian function for this problem is defined
as
\begin{equation}
{\cal L}_A ({\bf z},\blambda,\mu) = E({\bf z}) + \blambda^T ({\bf
Az-b}) + \frac{\mu}{2}\,  \|{\bf Az-b}\|_2^2,\label{augmented_L}
\end{equation}
where $\blambda \in \mathbb{R}^p$ is a vector of Lagrange
multipliers and $\mu \geq 0$ is called the AL penalty
parameter \cite{NocedalWright}. The AL algorithm iterates
between minimizing ${\cal L}_A ({\bf z},\blambda,\mu)$ with respect
to ${\bf z}$, keeping $\blambda$ fixed, and updating $\blambda$.

\vspace{0.1cm}
\begin{algorithm}{AL}{
\label{alg:salsa1}}
Set $k=0$, choose $\mu > 0$, ${\bf z}_0$, and  $\blambda_0$.\\
\qrepeat\\
     ${\bf z}_{k+1} \in \arg\min_{{\bf z}} {\cal L}_A ({\bf z},\blambda_k,\mu)$\\
     $\blambda_{k+1} \leftarrow \blambda_{k} + \mu ({\bf Az}_{k+1} - {\bf b})$\\
     $k \leftarrow k + 1$
\quntil stopping criterion is satisfied.
\end{algorithm}
\vspace{0.1cm}

It is possible (in some cases recommended) to update the
value of $\mu$ at each iteration \cite{NocedalWright}, \cite{Bazaraa}
(Chap.~9). However, unlike in the quadratic penalty method, it is
not necessary to take $\mu$ to infinity to guarantee that the AL
converges to the solution of the constrained problem
(\ref{constrained_linear}). In this paper, we will consider only
the case of fixed $\mu$.

After a straightforward manipulation, the terms
added to $E({\bf z})$ in  ${\cal L}_A ({\bf z},\blambda_k,\mu)$
(see (\ref{augmented_L})) can be written as a single quadratic term,
leading to the following alternative form for the AL algorithm:

\vspace{0.1cm}
\begin{algorithm}{AL (version 2)}{
\label{alg:salsa2}}
Set $k=0$, choose $\mu > 0$, ${\bf z}_0$, and  ${\bf d}_0$.\\
\qrepeat\\
     ${\bf z}_{k+1} \in \arg\min_{{\bf z}} E({\bf z}) + \frac{\mu}{2}\|{\bf Az-d}_k\|_2^2$\\
     ${\bf d}_{k+1} \leftarrow {\bf d}_{k} - ({\bf Az}_{k+1} - {\bf b})$\\
     $k \leftarrow k+1$
\quntil stopping criterion is satisfied.
\end{algorithm}
\vspace{0.1cm}

This form of the AL algorithm makes clear its equivalence with the
Bregman iterative method, as given in \cite{YinOsherGoldfarbDarbon}.

\subsection{AL for Variable Splitting and Its
Convergence} \label{ASAL}
Problem (\ref{constrained_basic}) can be written in the form
(\ref{constrained_linear})  with ${\bf z} =
[{\bf u}^{T}, \ {\bf v}^{T}]^T$, ${\bf b} = {\bf 0}$,
${\bf A} =  [\,{\bf I} \;\; -{\bf I}\,]$, and
$E({\bf z}) = f_1({\bf u}) + f_2({\bf v})$.
With these definitions in place, Steps 3 and 4 of the AL algorithm
(version 2) can be written as follows
\begin{equation}\label{mixed}
\left(\begin{array}{cc}
    {\bf u}_{k+1}\\
 {\bf v}_{k+1}
\end{array}\right) \in  \arg\min_{{\bf
u},{\bf v}} f_{1}({\bf u}) + f_{2}({\bf v}) + \frac{\mu}{2} \|{\bf u} - {\bf v} - {\bf d}_k\|_2^2
\end{equation}
\begin{equation}
{\bf d}_{k+1} = {\bf d}_{k} - ({\bf u}_{k+1} - {\bf v}_{k+1}).
\end{equation}
The minimization problem (\ref{mixed}) is clearly non-trivial: in
general, it involves non-separable quadratic and possibly non-smooth
terms. A natural approach is to use a non-linear
block-Gauss-Seidel (NLBGS) technique, in which (\ref{mixed}) is
solved by alternating  minimization with respect to ${\bf u}$ and
${\bf v}$, while keeping the other variable fixed. Remarkably, it
has been shown that the AL algorithm converges, even if the
exact solution of (\ref{mixed}) is replaced with a single NLBGS step
\cite[Theorem 8]{EcksteinBertsekas} (see also \cite{Setzer}).
The resulting algorithm  is as  follows.

\vspace{0.1cm}
\begin{algorithm}{Alternating Split AL}{
\label{alg:salsa2}}
Set $k=0$, choose $\mu > 0$, ${\bf u}_0$, ${\bf v}_0$, and  ${\bf d}_0$.\\
\qrepeat\\
   $  {\bf u}_{k+1}  \in  \arg\min_{{\bf u}} f_{1}({\bf u})
 + \frac{\mu}{2} \|{\bf u} - {\bf v}_k - {\bf d}_k\|_2^2$\\
  $  {\bf v}_{k+1}  \in  \arg\min_{{\bf v}} f_{2}({\bf v})
 + \frac{\mu}{2} \|{\bf u}_{k+1} - {\bf v} - {\bf d}_k\|_2^2$\\
     ${\bf d}_{k+1} \leftarrow {\bf d}_{k} - {\bf u}_{k+1} + {\bf v}_{k+1} $\\
     $k \leftarrow k+1$
\quntil stopping criterion is satisfied.
\end{algorithm}
\vspace{0.1cm}

Problem (\ref{problem}) has the form (\ref{unconstrained_basic})
where $f_1$ is quadratic, thus Step 3 consist in solving a linear
system of equations. We will return to the particular form of this system in
the next section. With $f_2 = \tau\phi$, a regularizer, Step 4
corresponds to applying a shrinkage/thresholding function, that is,
${\bf v}_{k+1}  =  \bPsi_{\tau\phi/\mu}\left({\bf u}_{k+1} - {\bf
d}_k\right),$ usually a computationally inexpensive operation.

\section{Proposed Method} \subsection{Constrained Optimization
Formulation and Algorithm} Performing the VS on
problem (\ref{problem}) yields the following constrained
formulation:
\begin{equation}\begin{array}{cl}
{\displaystyle \min_{\bbeta,\btheta}} & \frac{1}{2}\|{\bf HW} \bbeta - {\bf y}\|_2^2 + \tau \, \phi (\btheta) \\
\mbox{s.t.} & \bbeta = \btheta. \end{array}
\label{problem_constrained_image}
\end{equation}
This VS decouples the quadratic non-separable term
$\|{\bf HW} \bbeta - {\bf y}\|_2^2$ from the
non-quadratic term $\phi (\btheta)$, to deal with  the
non-separability of the quadratic data term. In contrast,
split-Bregman methods use a splitting to avoid non-separability of
the regularizer.

Problem (\ref{problem_constrained_image}) has the form
(\ref{constrained_basic}), with ${\bf u} = \bbeta$, ${\bf v}  = \btheta$,
$f_1({\bf u}) = (1/2)\|{\bf HW} \bbeta - {\bf y}\|_2^2$, and
$f_2({\bf v}) = \tau \, \phi (\btheta)$. Applying this translation
table to the {\it Alternating Split AL} algorithm presented Section
\ref{ASAL}, we obtain the following algorithm.

\vspace{0.1cm}
\begin{algorithm}{Split AL Shrinkage Algorithm}{
\label{alg:salsa2}}
Set $k=0$, choose $\mu > 0$, $\bbeta_0$, $\btheta_0$, and  ${\bf d}_0$.\\
\qrepeat\\
   $\bbeta'_{k} = \btheta_{k} + {\bf d}_k$ \\
   ${\displaystyle \bbeta_{k+1}  \in  \arg\min_{\bbeta}  \|{\bf HW} \bbeta - {\bf y}\|_2^2
 + \mu \|\bbeta - \bbeta'_k\|_2^2}$\\
 $\btheta'_k = \bbeta_{k+1} - {\bf d}_k$\\
  ${\displaystyle  \btheta_{k+1} = \bPsi_{\tau\phi/\mu}(\btheta'_k),}$\\
     ${\displaystyle {\bf d}_{k+1} \leftarrow {\bf d}_{k} - \bbeta_{k+1} + \btheta_{k+1}}$\\
     ${\displaystyle k \leftarrow k+1}$
\quntil stopping criterion is satisfied.
\end{algorithm}
\vspace{0.1cm}

Since Step 4 is a strictly convex quadratic problem, its
solution is unique and given by
\begin{equation}
\bbeta_{k+1} = \left({\bf W}^T{\bf H}^T{\bf H}{\bf W} + \mu\, {\bf
I}\right)^{-1}\left( {\bf W}^T{\bf H}^T{\bf y} + \mu\, \bbeta'_{k}
\right). \label{GS_linear}
\end{equation}
In the next subsection, we show how $\bbeta_{k+1}$ can be efficiently computed.
Note that $\left({\bf W}^T{\bf H}^T{\bf H}{\bf W} \! + \!\mu\, {\bf
I}\right)$ is a regularized (by the addition of $\mu {\bf I}$) version
of the Hessian of $\frac{1}{2}\|{\bf HW} \bbeta - {\bf y}\|_2^2$.

\subsection{Computing $\bbeta_{k+1}$}
Assume that ${\bf W}$ is a normalized tight (Parseval) frame,
{\it i.e.}, ${\bf W\,W}^T = {\bf I}$ (although
possibly  ${\bf W}^T{\bf W} \neq {\bf I}$), and that
${\bf H}$ is the matrix representation of a convolution,
{\it i.e.}, products by ${\bf H}$ or ${\bf H}^T$ can be computed
in the Fourier domain, with $O(n\log n)$ cost via the FFT.

The assumptions in the previous paragraph will enable us to
compute the matrix inversion in (\ref{GS_linear}), even if
it is not feasible to explicitly form matrix ${\bf HW}$.
Using the  Sherman–-Morrison-–Woodbury inversion formula,
(\ref{GS_linear}) becomes
\begin{eqnarray}
\bbeta_{k+1} & \!\!\! = \!\!\!  & \frac{1}{\mu} \left( {\bf I} -
{\bf W}^T {\bf H}^T \left({\bf H}{\bf W}{\bf W}^T {\bf H}^T + \mu\,
{\bf I}\right)^{-1} {\bf H} {\bf W}\right) {\bf r}_k
\nonumber \\
& \!\!\! = \!\!\!  & \frac{1}{\mu} \left( {\bf I} - {\bf W}^T {\bf
H}^T \left({\bf H} {\bf H}^T + \mu\, {\bf I}\right)^{-1} {\bf H}
{\bf W}\right) {\bf r}_k \label{GS_linear_3}
\end{eqnarray}
where ${\bf r}_k = \left( {\bf W}^T {\bf H}^T {\bf y} + \mu\,
\bbeta'_{k} \right)$. Furthermore, since ${\bf H}$ is the matrix
representation of a convolution, (\ref{GS_linear_3}) can be written as
\begin{equation}
\bbeta_{k+1} = \frac{1}{\mu} \left( {\bf I} - {\bf W}^T {\bf F}
{\bf W}\right) {\bf r}_k, \label{GS_linear_5}
\end{equation}
where
\begin{equation}
{\bf F} = {\bf U}^H {\bf D^*} \left( |{\bf D}|^2  + \mu\, {\bf
I}\right)^{-1} {\bf D U}, \label{matrix_fk}
\end{equation}
${\bf U}$ and ${\bf U}^H$ are the matrix representations of
the forward and inverse discrete Fourier transform (DFT), and
${\bf D}$ is a diagonal matrix containing the DFT of the
convolution represented by ${\bf H}$. Notice that
the product by ${\bf F}$ corresponds to applying a filter in
the DFT domain, which can be done using FFT algorithms with $O(n\log n)$
cost. Notice also that ${\bf H}^T{\bf W}^T {\bf y}$ can be precomputed.
When the products by ${\bf W}^T$ and ${\bf W}$ are direct
and inverse tight frame transforms for which fast algorithms exist,
the leading cost of each application of (\ref{GS_linear_5}) will be
either $O(n\log n)$ or the cost of these frame transforms (usually
also $O(n \log n)$).

Finally, the complete algorithm, which we term SALSA
(split augmented Lagrangian shrinkage algorithm) is as follows.

\vspace{0.1cm}
\begin{algorithm}{SALSA}{
\label{alg:salsa_complete}} {\bf Initialization:}  set $k=0$;
choose $\mu > 0$, $\bbeta_0$, $\btheta_0$,   ${\bf d}_0$;\\
compute $\bar{\bf y} = {\bf W}^T{\bf H}^T{\bf y}$\\
compute ${\bf F} = {\bf U}^H {\bf D^*} \left( |{\bf D}|^2  + \mu {\bf I}\right)^{-1} {\bf D U}$\\
\qrepeat\\
$\bbeta'_{k} \leftarrow \btheta_k + {\bf d}_k$\\
${\bf r}_k \leftarrow \bar{\bf y} + \mu \bbeta'_{k} $\\
$\bbeta_{k+1} \leftarrow \frac{1}{\mu} \left( {\bf I} - {\bf W}^T {\bf F}_k {\bf W}\right) {\bf r}_k$\\
$\theta'_k \leftarrow \bbeta_{k+1} - {\bf d}_k$\\
$\theta_{k+1} \leftarrow \bPsi_{\tau\phi/\mu} (\theta'_k)$\\
${\bf d}_{k+1} \leftarrow {\bf d}_{k} - \bbeta_{k+1} + \btheta_{k+1}$\\
$k \leftarrow k+1$ \quntil stopping  criterion is satisfied.
\end{algorithm}
\vspace{0.1cm}

\section{Experiments}
\begin{figure*}[ht]
\centering
\subfigure[]{
\includegraphics[width=0.3\textwidth]{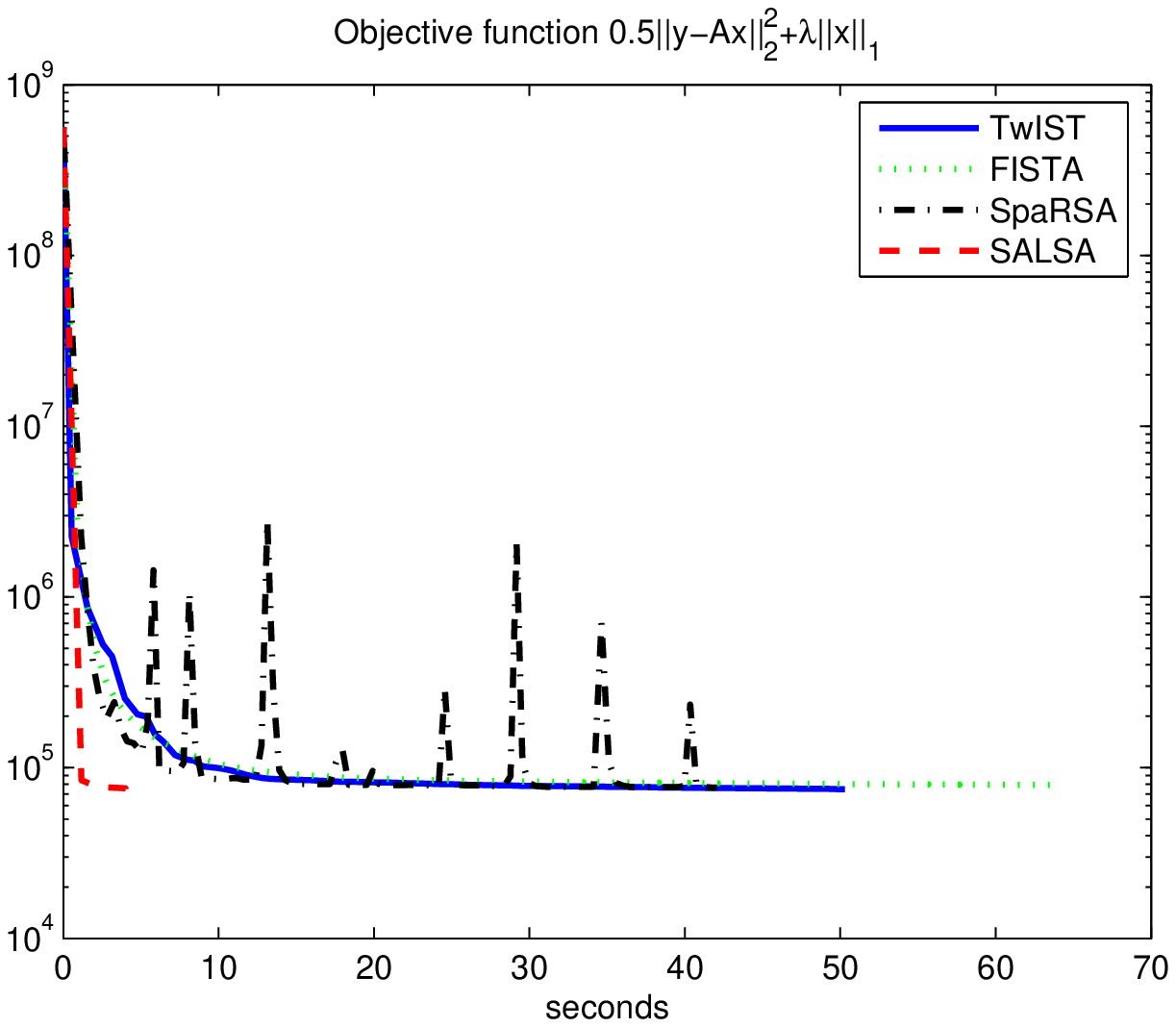}
}
\subfigure[]{
\includegraphics[width=0.3\textwidth]{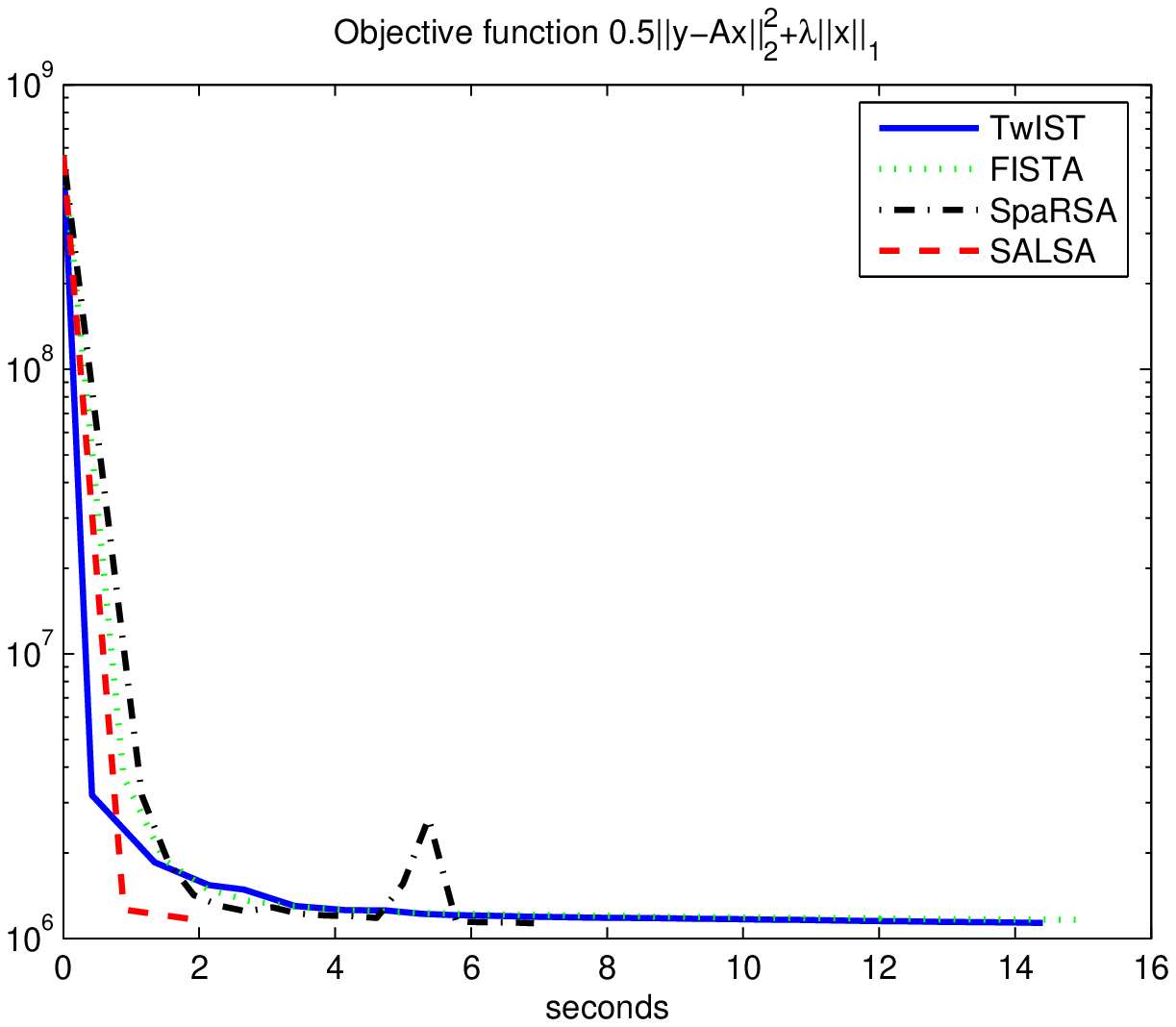}
}
\subfigure[]{
\includegraphics[width=0.3\textwidth]{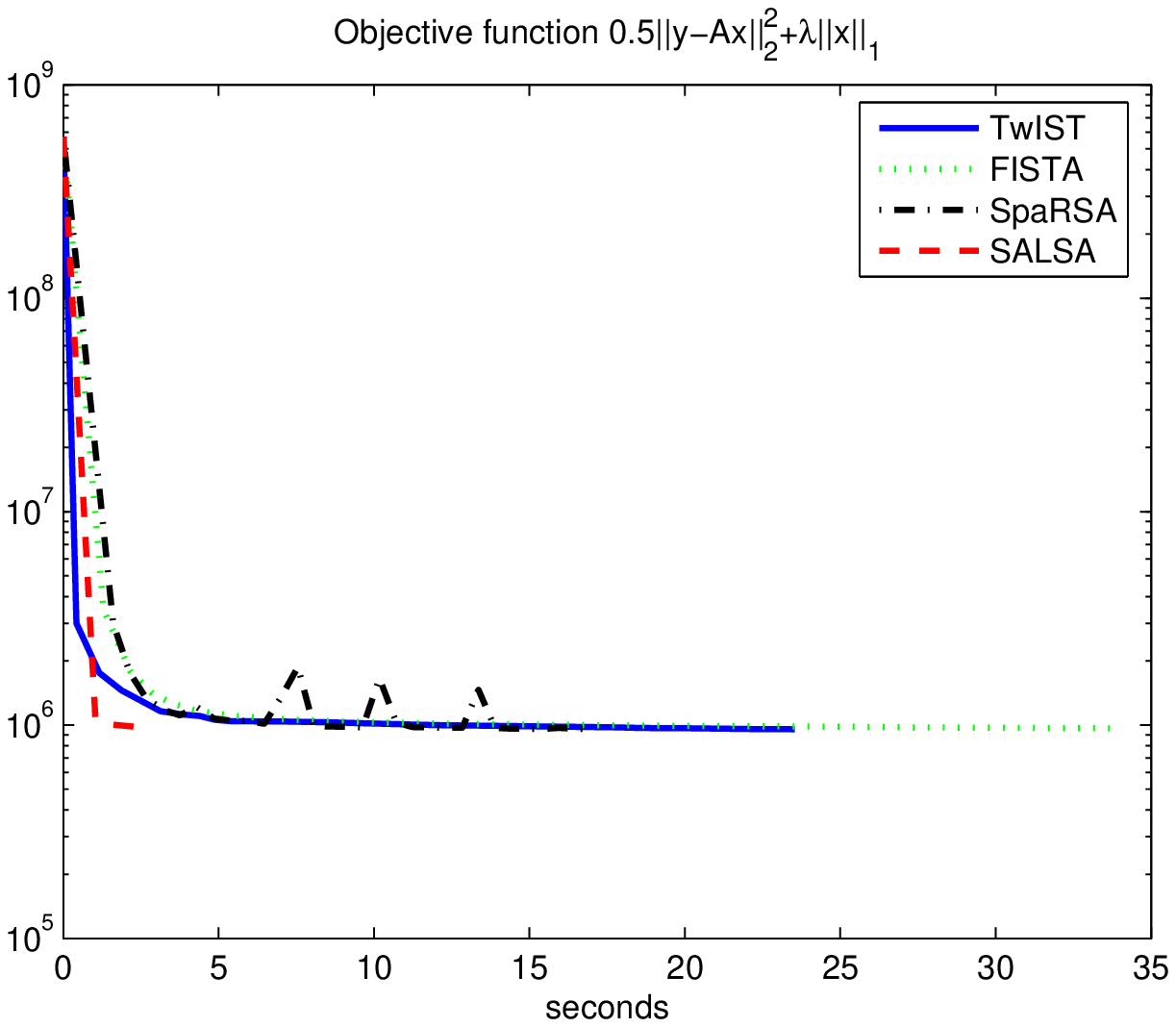}
}
\label{fig:evolutionobjective}
\caption{Objective function evolution: (a) $9\times 9$ uniform
blur, $\sigma=0.56$; (b) Gaussian blur, $\sigma^2=8$; (c) $h_{ij} = 1/(1 + i^2 + j^2)$ blur, $\sigma^2=2$.}
\end{figure*}

We consider five standard image deconvolution benchmark problems \cite{FigueiredoNowak2003},
summarized in Table~\ref{decon_problems}, all on the well-known Cameraman image. The regularizer
is $\phi(\bbeta) = \|\bbeta\|_1$, thus $\bPsi_{\tau\phi}$ is a soft threshold.
In all the experiments, ${\bf W}$ is a redundant Haar wavelet frame, with $4$ levels,
and the blur operator ${\bf H}$ is applied via the FFT. The regularization parameter
$\tau$ in each case was hand tuned for best improvement in SNR.
The value of $\mu$ for fastest convergence was found to differ in each case, but a
good rule of thumb, used in all the experiments, is $\mu = 0.1\tau$.
We compare SALSA with current state of the art methods: TwIST \cite{TwIST},
SpaRSA \cite{SpaRSA_SP}, and FISTA \cite{FISTA},
in terms of the time taken to reach the same value of the objective function.
Table~\ref{tab:deconresultsl1redundant} shows the CPU times taken by each of
the algorithms in each of the experiments. Figure~\ref{fig:evolutionobjective}
shows the plots of the objective function
$\frac{1}{2}\|{\bf HW\beta - y}\|+\tau\|{\bf \beta}\|_1$, evolving over time,
in experiments $1$, $2$B, and $3$A.

\begin{table}[h]
\centering
\caption{Details of the image deconvolution experiments.}\label{decon_problems}
\begin{tabular}{|c|l|l|}
  \hline
 Experiment &  blur kernel  \rule[-0.1cm]{0cm}{0.4cm}  & $\sigma^2$ \\ \hline
1  &$9\times 9$ uniform & $0.56^2$ \\
2A & Gaussian & 2\\
2B & Gaussian & 8\\
3A  & $h_{ij} = 1/(1 + i^2 + j^2)$ & 2 \\
3B  & $h_{ij} = 1/(1 + i^2 + j^2)$ & 8 \\
\hline
\end{tabular}
\end{table}

\begin{table}[h]
\centering \caption{CPU times (in seconds) for the various algorithms.}
\label{tab:deconresultsl1redundant}
\begin{tabular}{|c|l|l|l|l|}
  \hline
 Experiment & TwIST   \rule[-0.1cm]{0cm}{0.4cm}  & SpARSA  & FISTA & SALSA \\
\hline
1 & 50.2969 & 42.0469 & 64.2344 & 4.000 \\
2A & 30.7656 & 40.6094 & 61.7031 & 4.03125 \\
2B & 14.4063 & 6.92188 & 15.0781 & 1.9375 \\
3A & 23.5313 & 17.0156 & 33.7969 & 2.60938 \\
3B & 8.1875 & 6.17188 & 18.0781 & 1.89063 \\
\hline
\end{tabular}
\end{table}

\section{Conclusions} \label{sec:conclusions}
We have proposed a fast algorithm for frame-based image deconvolution,
based on variable splitting and solving the constrained optimization
problem through an augmented Lagrangian scheme. Experimental results
with $\ell_1$ regularization show that our new algorithm outperforms
existing state-of-the-art methods in terms of computation time,
by a considerable factor. Future work includes the application of
SALSA to other inverse problems, namely compressed sensing and
reconstruction with missing samples.

\end{document}